\documentclass{article}
\usepackage{graphicx} 

\usepackage[a4paper,top=2cm,bottom=2cm,left=3cm,right=3cm,marginparwidth=1.75cm]{geometry}

\usepackage[normalem]{ulem}

\usepackage[utf8]{inputenc} 
\usepackage[T2A]{fontenc}
\usepackage[english]{babel}
\usepackage[unicode, pdftex]{hyperref}

\usepackage{amsfonts}
\usepackage{amsthm}

\usepackage{amsmath}
\usepackage{graphicx}
\usepackage{tikz}
\usetikzlibrary{arrows}
\usetikzlibrary{arrows.meta}
\usepackage{caption}
\usepackage{subcaption}
\usepackage{comment}

\usepackage{algorithm}
\usepackage{algpseudocode}

\usetikzlibrary{arrows,automata,positioning}

\DeclareRobustCommand{\rchi}{{\mathpalette\irchi\relax}}
\newcommand{\irchi}[2]{\raisebox{\depth}{$#1\chi$}} 

\title{Decomposition of surface into quadrilaterals}
\date{2023\\ April}
\author{Kostiantyn Cherkashyn,\\  Tarash Shevchenko National University of Kyiv}

\begin{document}

\maketitle

\begin{abstract}
There are many approaches to the classification of Morse functions and their gradient fields (Morse Fields) on 2-surfaces. This paper studies the gluings of quadrilaterals and the classification of topological surfaces obtained by gluing together sides of quadrilateral graphs created by trajectories on gradient fields with focal critical points as vertices and a saddle point in center for classification saddle and focus critical points. We review the graphs embedded into surfaces associated with gluings, the number of possible gluings, and some algorithms based on labeling schemes of fundamental polygons that can be used to compute those numbers. 
\end{abstract}

\section{Introduction}

The topic of this paper refers to applications of topological graph theory to the classification of Morse flows on surfaces. It uses both topological and algebraic methods when working with embeddings of graphs on the surface.
The main construction, which dates back to the works of Peixoto, consists in the construction of a discriminant graph, the vertices of which are sources, and the edges of which are one-dimensional flowing manifolds.

Classic papers on the topological theory of graphs are books
\cite{angelini2018linear, appel1976every, bang2008digraphs, bollobas1998modern, bondy2008graph, chrobak1994planar, diestel2005graph, dillencourt1993graph, eppstein1999quasi, godsil2001algebraic, gross1999topological, liotta2018beyond, lovasz2001graph, mohar1991embeddings, mohar2001graphs, neumann2015topological, pach2018planar, thomassen2009graphs}.

Embedded  graphs as topological invariants of flows on closed surfaces were constructed in  \cite{bilun2023gradient, Kybalko2018, Oshemkov1998, Peixoto1973, prishlyak1997graphs, prishlyak2020three, akchurin2022three, prishlyak2022topological, prishlyak2017morse, kkp2013, prishlyak2021flows, prishlyak2020topology, prishlyak2019optimal, prishlyak2022Boy},
and on surfaces with a boundary y
\cite{bilun2023discrete, bilun2023typical, loseva2016topology, loseva2022topological, prishlyak2017morse, prishlyak2022topological, prishlyak2003sum, prishlyak2003topological, prishlyak1997graphs, prishlyak2019optimal, stas2023structures}.
For a 3-manifold, a Hegaard diagram is an embedded 4-valent graph in a surface \cite{prish1998vek, prish2001top, Prishlyak2002beh2, prishlyak2002ms, prishlyak2007complete, hatamian2020heegaard, bilun2022morse, bilun2022visualization}.

Morse streams are gradient streams for Morse functions. If we fix the value of the function at special points, then the structure of the flow determines the structure of the function \cite{lychak2009morse, Smale1961}.

Topological invariants, as Reeb graphs, of functions on oriented surfaces were constructed in \cite{Kronrod1950} and \cite{Reeb1946} and in \cite{lychak2009morse} for unoriented surfaces, and in \cite{Bolsinov2004, hladysh2017topology, hladysh2019simple, prishlyak2012topological} for of surfaces with a boundary, in \cite{prishlyak2002morse} for non-compact surfaces.

Embedded graphs as topological invariants of smooth functions were also studied in papers \cite{bilun2023morseRP2, bilun2023morse, hladysh2019simple, hladysh2017topology, prishlyak2002morse, prishlyak2000conjugacy, prishlyak2007classification, lychak2009morse, prishlyak2002ms, prish2015top, prish1998sopr, bilun2002closed, Sharko1993}, for manifolds with a boundary in papers \cite {hladysh2016functions, hladysh2019simple, hladysh2020deformations}, and on 3- and 4-dimensional manifolds in \cite{prishlyak1999equivalence, prishlyak2001conjugacy}.

To get acquainted with the topological theory of functions and dynamical systems, we recommend  \cite{prishlyak2012topological, prish2002theory, prish2004difgeom, prish2006osnovy, prish2015top}.
\par
The purpose of this paper is to study topological surfaces created by gluing together sides of quadrilaterals, a variation of the classical question that arises in many different applications: in how many ways edges of 2N-gon can be pairwise glued to obtain a certain surface? We consider complete and incomplete i.e. some pairs of edges remain unglued, orientable, and non-orientable, but only connected surfaces. The main focus of this research is algorithms acting on labeling schemes, classification, and enumeration of labeling schemes equivalence classes under actions of the dihedral group $D_n$.
\par
General solutions to a similar problem can be found in \cite{Harer1986euler} for closed orientable surfaces, considering rotations and reflections of the fundamental polygon as a different gluing. And in \cite{akhmedov2008gluing} there is a solution for orientable surfaces including surfaces with one or more boundary components, considering surface rotation as a single equivalence class, but preserving orientation in the sense that reflection symmetry of the polygon is considered as a different embedded graph.
\par

\section{Definitions}

Let $X$ be a polygonal region, the \textit{gluing} of two distinct edges $a$ and $b$ of $X$ is a quotient space of this polygon $$X/\sim, \ g(x) \sim h(x) \ or \ g(x) \sim h(1-x), \quad x \in [0,1]$$ where 

$g: [0,1] \rightarrow X$ is a path in $X$ that consists of all points of the edge $a$ and 
$h: [0,1] \rightarrow X$ is a path corresponds to edge $b$.

\begin{figure}[!ht]
\centering
\begin{tikzpicture}[]
\node (A) at (0,2) {\textbullet};
\node (B) at (1,2) {\textbullet};
\node (C) at (1.75,1) {\textbullet};
\node (D) at (1,0) {\textbullet};
\node (E) at (0,0) {\textbullet};
\node (F) at (-0.75,1) {\textbullet};

\node at (-0.2,2.2) {$a_1$};
\node at (1.2,2.2) {$a_2$};
\node at (-0.2,-0.2) {$b_1$};
\node at (1.2,-0.2) {$b_2$};

\draw [-](A) edge node[above] {a} (B);
\draw [-](B) edge node[above right] {} (C);
\draw [-](C) edge node[below right] {} (D);
\draw [-](D) edge node[below] {b} (E);
\draw [-](E) edge node[below left] {} (F);
\draw [-](F) edge node[above left] {} (A);

\node at(2.5,1.1) {$\xrightarrow{f}$};

\node (Z) at (3.5,1) {\textbullet};
\node (X) at (4.5,1) {\textbullet};
\node (B) at (6,1) {\textbullet};
\node (N) at (7,1) {\textbullet};

\draw [-](Z) to [bend right] node[above] {} (X);
\draw [-](Z) to [bend left] node[above] {} (X);
\draw [-](X) to  node[above] {$f(b)$} (B);
\draw [-](B) to [bend right] node[above] {} (N);
\draw [-](B) to [bend left] node[above] {} (N);

\node (X) at (4.5,0.5) {$f(b_1)$};
\node (B) at (6,0.5) {$f(b_2)$};

\end{tikzpicture}
    \caption{Graph amalgamation}
    \label{fig:amalgamation}
\end{figure}
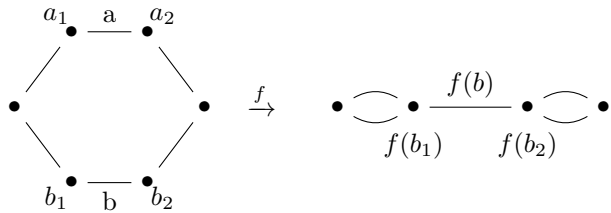

 Two equivalence relations $g(x) \sim h(x)$ and $g(x) \sim h(1-x)$ correspond to two possible ways to glue two edges that may result in orientable or non-orientable. In respect to graph $G$ that represents a polygonal boundary, this gluing is the amalgamation of the graph $G_f$, where $f$ is a graph isomorphism mapping edges $a \leftrightarrow b$ and their vertices $a_1 \leftrightarrow b_1$ and $a_2 \leftrightarrow b_2$ or $a_1 \leftrightarrow b_2$ and $a_2 \leftrightarrow b_1$ \ref{fig:amalgamation}. In this paper, we calculate the number of different embedded graphs created this way on corresponding surfaces

\par
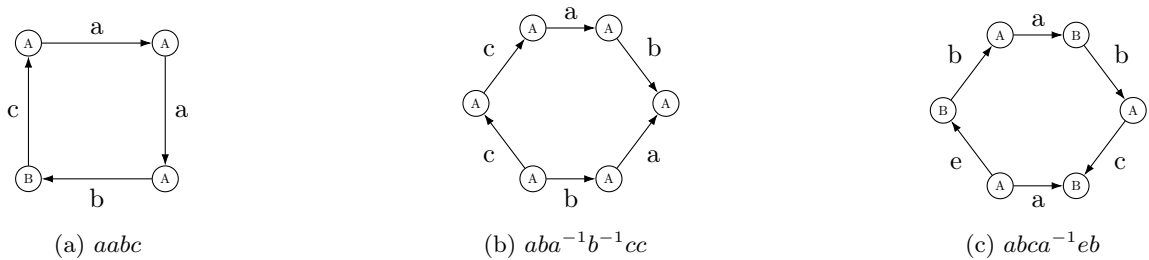
\begin{figure}[!ht]
    \begin{subfigure}{.18\linewidth}
\centering
\begin{tikzpicture}[>={Latex[width=1mm,length=1.5mm]},
                               node distance = 1.5cm and 2cm,
                               el/.style = {inner sep=2pt, align=center, sloped},
                               every label/.append style = {font=\tiny}]
\node[shape=circle,draw=black,scale=0.5] (A) at (0,1.8) {A};
\node[shape=circle,draw=black,scale=0.5] (B) at (1.8,1.8) {A};
\node[shape=circle,draw=black,scale=0.5] (C) at (1.8,0) {A};
\node[shape=circle,draw=black,scale=0.5] (D) at (0,0) {B};

\path [->](A) edge node[above] {a} (B);
\path [->](B) edge node[right] {a} (C);
\path [->](C) edge node[below] {b} (D);
\path [->](D) edge node[left] {c} (A);
\end{tikzpicture}
\caption{$a a b c $}
\end{subfigure}
 \hfill 
 \begin{subfigure}{.18\linewidth}
\begin{tikzpicture}[>={Latex[width=1mm,length=1.5mm]},
                                node distance = 3cm and 4cm,
                                el/.style = {inner sep=2pt, align=center, sloped},
                                every label/.append style = {font=\tiny}]
\node[shape=circle,draw=black,scale=0.5] (A) at (0.0,2.0) {A};
\node[shape=circle,draw=black,scale=0.5] (B) at (1.0,2.0) {A};
\node[shape=circle,draw=black,scale=0.5] (C) at (1.75,1.0) {A};
\node[shape=circle,draw=black,scale=0.5] (D) at (1.0,0.0) {A};
\node[shape=circle,draw=black,scale=0.5] (E) at (0.0,0.0) {A};
\node[shape=circle,draw=black,scale=0.5] (F) at (-0.75,1.0) {A};

\path [->](A) edge node[above] {a} (B);
\path [->](B) edge node[above right] {b} (C);
\path [<-](C) edge node[below right] {a} (D);
\path [<-](D) edge node[below] {b} (E);
\path [->](E) edge node[below left] {c} (F);
\path [->](F) edge node[above left] {c} (A);
\end{tikzpicture}
\caption{$a b a^{-1} b^{-1} c c $}
\end{subfigure}
\hfill
 \begin{subfigure}{.18\linewidth}
\begin{tikzpicture}[>={Latex[width=1mm,length=1.5mm]},
                                node distance = 3cm and 4cm,
                                el/.style = {inner sep=2pt, align=center, sloped},
                                every label/.append style = {font=\tiny}]
\node[shape=circle,draw=black,scale=0.5] (A) at (0.0,2.0) {A};
\node[shape=circle,draw=black,scale=0.5] (B) at (1.0,2.0) {B};
\node[shape=circle,draw=black,scale=0.5] (C) at (1.75,1.0) {A};
\node[shape=circle,draw=black,scale=0.5] (D) at (1.0,0.0) {B};
\node[shape=circle,draw=black,scale=0.5] (E) at (0.0,0.0) {A};
\node[shape=circle,draw=black,scale=0.5] (F) at (-0.75,1.0) {B};

\path [->](A) edge node[above] {a} (B);
\path [->](B) edge node[above right] {b} (C);
\path [->](C) edge node[below right] {c} (D);
\path [<-](D) edge node[below] {a} (E);
\path [->](E) edge node[below left] {e} (F);
\path [->](F) edge node[above left] {b} (A);
\end{tikzpicture}
\caption{$a b c a^{-1} e b $}
\end{subfigure}
    \caption{Gluing schemes}
    \label{fig:gluing_schemes}
\end{figure}

A common and useful representation of such gluings is \textit{gluing schemes}: a polygon with the labeled edges and orientation marked on each edge, an edge with the same label appears no more than twice\ref{fig:gluing_schemes}. For the purpose of this paper, we consider all rotations and reflection symmetries of the polygon to be equivalent i.e. invariant under all actions of the group $D_{2n}$ for \textit{2n}-gon.

\par
Another representation of gluing is a \textit{labeling scheme}: a word along the boundary of the polygon consequently labeling all the edges, if some edge $a$ is oriented against the chosen direction it is written as an inverse $a^{-1}$ if the letter and it's inverse appears in a single word it corresponds to non-orientable gluing of the corresponding edges. Letter and/or its inverse cannot appear in the word more than twice, letters that appear only once correspond to the free edges i.e. boundary of the surface.
\par
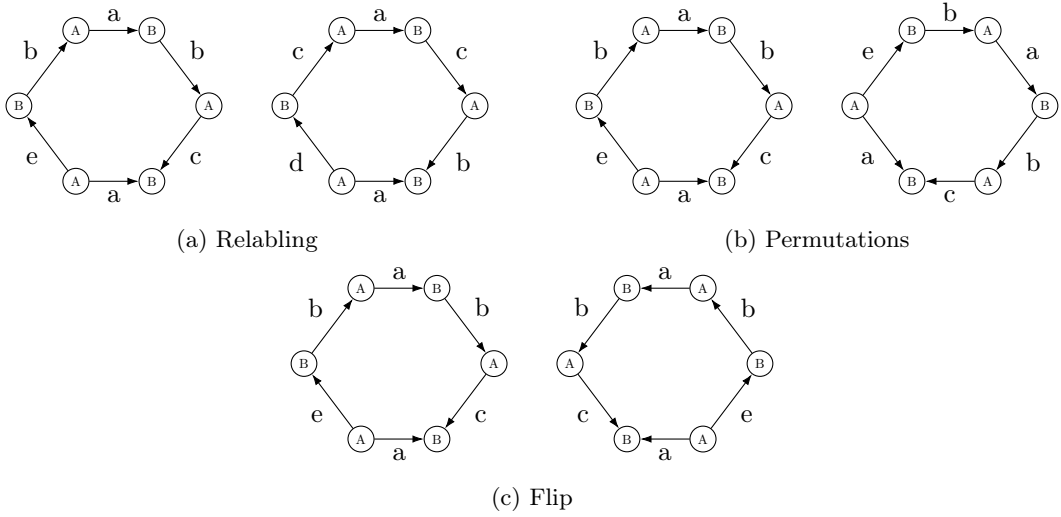
\begin{figure}[!ht]
    \centering
    \begin{subfigure}{.5\textwidth}
    \centering
        \begin{tikzpicture}[>={Latex[width=1mm,length=1.5mm]},
                               node distance = 1.5cm and 2cm,
                               el/.style = {inner sep=1pt, align=center, sloped},
                               every label/.append style = {font=\tiny}]
\node[shape=circle,draw=black, scale=0.5] (A) at (0,2) {A};
\node[shape=circle,draw=black, scale=0.5] (B) at (1,2) {B};
\node[shape=circle,draw=black, scale=0.5] (C) at (1.75,1) {A};
\node[shape=circle,draw=black, scale=0.5] (D) at (1,0) {B};
\node[shape=circle,draw=black, scale=0.5] (E) at (0,0) {A};
\node[shape=circle,draw=black, scale=0.5] (F) at (-0.75,1) {B};

\path [->](A) edge node[above] {a} (B);
\path [->](B) edge node[above right] {b} (C);
\path [->](C) edge node[below right] {c} (D);
\path [<-](D) edge node[below] {a} (E);
\path [->](E) edge node[below left] {e} (F);
\path [->](F) edge node[above left] {b} (A);

\node[shape=circle,draw=black, scale=0.5] (A) at (3.5,2) {A};
\node[shape=circle,draw=black, scale=0.5] (B) at (4.5,2) {B};
\node[shape=circle,draw=black, scale=0.5] (C) at (5.25,1) {A};
\node[shape=circle,draw=black, scale=0.5] (D) at (4.5,0) {B};
\node[shape=circle,draw=black, scale=0.5] (E) at (3.5,0) {A};
\node[shape=circle,draw=black, scale=0.5] (F) at (2.75,1) {B};

\path [->](A) edge node[above] {a} (B);
\path [->](B) edge node[above right] {c} (C);
\path [->](C) edge node[below right] {b} (D);
\path [<-](D) edge node[below] {a} (E);
\path [->](E) edge node[below left] {d} (F);
\path [->](F) edge node[above left] {c} (A);

\end{tikzpicture}
        \caption{Relabling}
        \label{fig:relabel}
    \end{subfigure}%
    \begin{subfigure}{.5\textwidth}
    \centering
        \begin{tikzpicture}[>={Latex[width=1mm,length=1.5mm]},
                               node distance = 1.5cm and 2cm,
                               el/.style = {inner sep=1pt, align=center, sloped},
                               every label/.append style = {font=\tiny}]
\node[shape=circle,draw=black, scale=0.5] (A) at (0,2) {A};
\node[shape=circle,draw=black, scale=0.5] (B) at (1,2) {B};
\node[shape=circle,draw=black, scale=0.5] (C) at (1.75,1) {A};
\node[shape=circle,draw=black, scale=0.5] (D) at (1,0) {B};
\node[shape=circle,draw=black, scale=0.5] (E) at (0,0) {A};
\node[shape=circle,draw=black, scale=0.5] (F) at (-0.75,1) {B};

\path [->](A) edge node[above] {a} (B);
\path [->](B) edge node[above right] {b} (C);
\path [->](C) edge node[below right] {c} (D);
\path [<-](D) edge node[below] {a} (E);
\path [->](E) edge node[below left] {e} (F);
\path [->](F) edge node[above left] {b} (A);

\node[shape=circle,draw=black, scale=0.5] (A) at (3.5,2) {B};
\node[shape=circle,draw=black, scale=0.5] (B) at (4.5,2) {A};
\node[shape=circle,draw=black, scale=0.5] (C) at (5.25,1) {B};
\node[shape=circle,draw=black, scale=0.5] (D) at (4.5,0) {A};
\node[shape=circle,draw=black, scale=0.5] (E) at (3.5,0) {B};
\node[shape=circle,draw=black, scale=0.5] (F) at (2.75,1) {A};

\path [->](A) edge node[above] {b} (B);
\path [->](B) edge node[above right] {a} (C);
\path [->](C) edge node[below right] {b} (D);
\path [->](D) edge node[below] {c} (E);
\path [<-](E) edge node[below left] {a} (F);
\path [->](F) edge node[above left] {e} (A);

\end{tikzpicture}
        \caption{Permutations}
        \label{fig:permutations}
    \end{subfigure}
    \begin{subfigure}{.5\textwidth}
    \centering
        \begin{tikzpicture}[>={Latex[width=1mm,length=1.5mm]},
                               node distance = 1.5cm and 2cm,
                               el/.style = {inner sep=1pt, align=center, sloped},
                               every label/.append style = {font=\tiny}]
\node[shape=circle,draw=black, scale=0.5] (A) at (0,2) {A};
\node[shape=circle,draw=black, scale=0.5] (B) at (1,2) {B};
\node[shape=circle,draw=black, scale=0.5] (C) at (1.75,1) {A};
\node[shape=circle,draw=black, scale=0.5] (D) at (1,0) {B};
\node[shape=circle,draw=black, scale=0.5] (E) at (0,0) {A};
\node[shape=circle,draw=black, scale=0.5] (F) at (-0.75,1) {B};

\path [->](A) edge node[above] {a} (B);
\path [->](B) edge node[above right] {b} (C);
\path [->](C) edge node[below right] {c} (D);
\path [<-](D) edge node[below] {a} (E);
\path [->](E) edge node[below left] {e} (F);
\path [->](F) edge node[above left] {b} (A);

\node[shape=circle,draw=black, scale=0.5] (A) at (3.5,2) {B};
\node[shape=circle,draw=black, scale=0.5] (B) at (4.5,2) {A};
\node[shape=circle,draw=black, scale=0.5] (C) at (5.25,1) {B};
\node[shape=circle,draw=black, scale=0.5] (D) at (4.5,0) {A};
\node[shape=circle,draw=black, scale=0.5] (E) at (3.5,0) {B};
\node[shape=circle,draw=black, scale=0.5] (F) at (2.75,1) {A};

\path [<-](A) edge node[above] {a} (B);
\path [<-](B) edge node[above right] {b} (C);
\path [<-](C) edge node[below right] {e} (D);
\path [->](D) edge node[below] {a} (E);
\path [<-](E) edge node[below left] {c} (F);
\path [<-](F) edge node[above left] {b} (A);

\end{tikzpicture}
        \caption{Flip}
        \label{fig:flip}
    \end{subfigure}
\end{figure}

Labeling schemes have some operations defined on them that do not affect it's resulting surface. In this paper we are going to use only a subset of them, namely:
\begin{enumerate}
\item \textit{Relabel}: change every occurrence of the same letter in the word to another letter, keeping exponents as it is. It is possible to flip all exponents of a given letter e.g. if $a$ is changed to $a^{-1}$ if there was an occurrence of $a^{-1}$ it has to be changed to $a$ \ref{fig:relabel}.
\item \textit{Permute}: cyclic permutation of all letters in the word e.g. $aba^{-1}b^{-1}$ can be changed to $ba^{-1}b^{-1}a$ or $a^{-1}b^{-1}ab$. This operation corresponds to a rotation of a polygon if we fix in space the location of the first edge \ref{fig:permutations}.
\item \textit{Flip}: reverse the word and flip all exponents e.g. $aba^{-1}b^{-1}$ is flipped to $bab^{-1}a^{-1}$. This operation corresponds to flipping over the polygon i.e. reflection symmetry \ref{fig:flip}.
\end{enumerate}

\newpage
\section{Classification of Surfaces}
Classification of surfaces is a very important and well-studied topic in topology. For closed surfaces, it is summarized by classification theorem.
\par
Let $S$ be a surface. $S$ is homeomorphic to precisely one of the following surfaces:
    $$S^2 \# \underbrace{T \# T \#...\# T}_n, \quad n \ge 0$$
    $$S^2 \# \underbrace{\mathbb{R}P^2 \# \mathbb{R}P^2 \#...\# \mathbb{R}P^2}_m, \quad m > 0$$

Proof of this theorem to some degree can be found in many different textbooks e.g. \cite{Kosniowski1980first}
\par
The above surfaces can be represented with standard labeling schemes:
$$ S^2 \ \simeq \ aa^{-1} $$
$$ \underbrace{T \# T \# ... \# T}_n \ \simeq \ a_1b_1a^{-1}_1b^{-1}_1a_2b_2a^{-1}_2b^{-1}_2...a_nb_na^{-1}_nb^{-1}_n $$
$$\underbrace{\mathbb{R}P^2 \# \mathbb{R}P^2 \#...\# \mathbb{R}P^2}_m \ \simeq \ a_1a_1a_2a_2...a_ma_m$$
\par
To complete this definition for surfaces with boundaries we introduce labeling schemes of surfaces containing a single occurrence of some letters without the occurrence of an inverse, which corresponds to surfaces with some number of discs removed i.e. some number of boundary components.
\par
We classify the surfaces based on the following topological invariants:
\begin{enumerate}
    \item Euler characteristic $\rchi = V - E + F$, where $E$ is the number of edges of the embedded graph, $V$ is the number of vertices, $F$ is the number of faces. It can also be expressed as $\rchi = 2 - 2G - B$ for orientable surface and $\rchi = 2 - G - B$ for non-orientable surface where $G$ is the genus for orientable surface or the demigenus for non-orientable surface and $B$ is the number of boundary components.
    \item Orientability: a topological property of a surface that can be vaguely defined as a surface having two sides (orientable) or only one side (non-orientable). A surface is non-orientable if it contains a subset homeomorphic to a Möbius strip i.e. has a cross-cap. The labeling scheme is of non-orientable type if it contains two occurrences of the same letter with the same exponent. If at least one such occurrence exists, every other glued pair except a cancelable type $aa^{-1}$ adds an additional cross-cap to a surface decreasing its Euler characteristic by one.
    \item Number of boundary components: a surface with $n$ boundary components is homeomorphic to a surface with $n$ discs removed.
\end{enumerate}
\par
Together listed above topological invariants can be used to determine a topological type of a surface and can be easily derived from the labeling schemes.

\section{Embedded Graphs}

In this paper, we are not only classifying surfaces obtained by gluing sides of polygons but calculating a number of ways to glue these surfaces for some number of quadrilaterals. For connected surfaces gluings from $n$ quadrilaterals are equivalent to the gluings from $2n+2$ sided polygons, but these gluings are introducing additional edges on the embedded graphs originating from boundaries of initial quadrilaterals.
\par
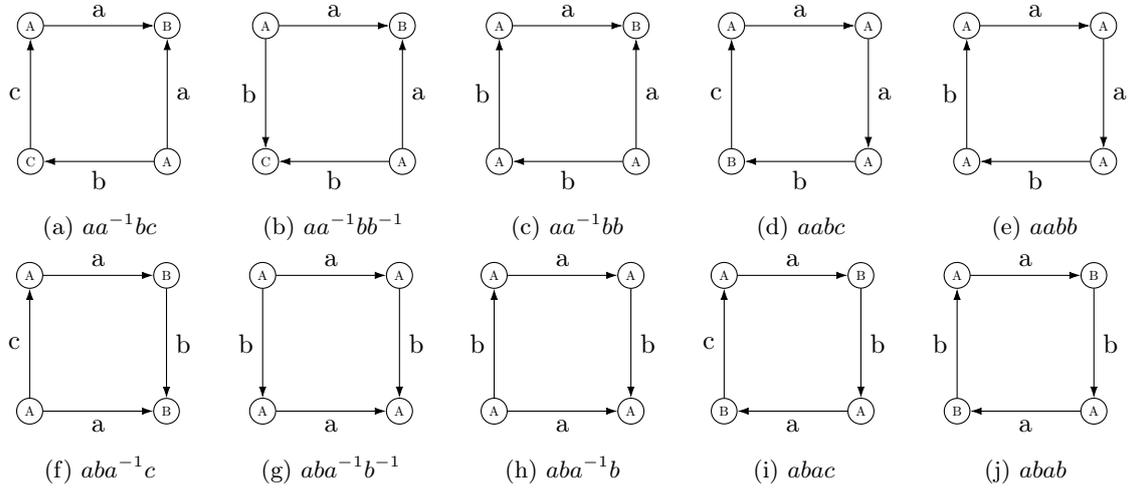
\begin{figure}[!ht]
    \begin{subfigure}{.18\linewidth}
\centering
\begin{tikzpicture}[>={Latex[width=1mm,length=1.5mm]},
                               node distance = 1.5cm and 2cm,
                               el/.style = {inner sep=2pt, align=center, sloped},
                               every label/.append style = {font=\tiny}]
\node[shape=circle,draw=black,scale=0.5] (A) at (0,1.8) {A};
\node[shape=circle,draw=black,scale=0.5] (B) at (1.8,1.8) {B};
\node[shape=circle,draw=black,scale=0.5] (C) at (1.8,0) {A};
\node[shape=circle,draw=black,scale=0.5] (D) at (0,0) {C};

\path [->](A) edge node[above] {a} (B);
\path [<-](B) edge node[right] {a} (C);
\path [->](C) edge node[below] {b} (D);
\path [->](D) edge node[left] {c} (A);
\end{tikzpicture}
\caption{$a a^{-1} b c $}
\end{subfigure}
 \hfill 
\begin{subfigure}{.18\linewidth}
\centering
\begin{tikzpicture}[>={Latex[width=1mm,length=1.5mm]},
                               node distance = 1.5cm and 2cm,
                               el/.style = {inner sep=2pt, align=center, sloped},
                               every label/.append style = {font=\tiny}]
\node[shape=circle,draw=black,scale=0.5] (A) at (0,1.8) {A};
\node[shape=circle,draw=black,scale=0.5] (B) at (1.8,1.8) {B};
\node[shape=circle,draw=black,scale=0.5] (C) at (1.8,0) {A};
\node[shape=circle,draw=black,scale=0.5] (D) at (0,0) {C};

\path [->](A) edge node[above] {a} (B);
\path [<-](B) edge node[right] {a} (C);
\path [->](C) edge node[below] {b} (D);
\path [<-](D) edge node[left] {b} (A);
\end{tikzpicture}
\caption{$a a^{-1} b b^{-1} $}
\end{subfigure}
 \hfill 
\begin{subfigure}{.18\linewidth}
\centering
\begin{tikzpicture}[>={Latex[width=1mm,length=1.5mm]},
                               node distance = 1.5cm and 2cm,
                               el/.style = {inner sep=2pt, align=center, sloped},
                               every label/.append style = {font=\tiny}]
\node[shape=circle,draw=black,scale=0.5] (A) at (0,1.8) {A};
\node[shape=circle,draw=black,scale=0.5] (B) at (1.8,1.8) {B};
\node[shape=circle,draw=black,scale=0.5] (C) at (1.8,0) {A};
\node[shape=circle,draw=black,scale=0.5] (D) at (0,0) {A};

\path [->](A) edge node[above] {a} (B);
\path [<-](B) edge node[right] {a} (C);
\path [->](C) edge node[below] {b} (D);
\path [->](D) edge node[left] {b} (A);
\end{tikzpicture}
\caption{$a a^{-1} b b $}
\end{subfigure}
 \hfill 
\begin{subfigure}{.18\linewidth}
\centering
\begin{tikzpicture}[>={Latex[width=1mm,length=1.5mm]},
                               node distance = 1.5cm and 2cm,
                               el/.style = {inner sep=2pt, align=center, sloped},
                               every label/.append style = {font=\tiny}]
\node[shape=circle,draw=black,scale=0.5] (A) at (0,1.8) {A};
\node[shape=circle,draw=black,scale=0.5] (B) at (1.8,1.8) {A};
\node[shape=circle,draw=black,scale=0.5] (C) at (1.8,0) {A};
\node[shape=circle,draw=black,scale=0.5] (D) at (0,0) {B};

\path [->](A) edge node[above] {a} (B);
\path [->](B) edge node[right] {a} (C);
\path [->](C) edge node[below] {b} (D);
\path [->](D) edge node[left] {c} (A);
\end{tikzpicture}
\caption{$a a b c $}
\end{subfigure}
 \hfill 
\begin{subfigure}{.18\linewidth}
\centering
\begin{tikzpicture}[>={Latex[width=1mm,length=1.5mm]},
                               node distance = 1.5cm and 2cm,
                               el/.style = {inner sep=2pt, align=center, sloped},
                               every label/.append style = {font=\tiny}]
\node[shape=circle,draw=black,scale=0.5] (A) at (0,1.8) {A};
\node[shape=circle,draw=black,scale=0.5] (B) at (1.8,1.8) {A};
\node[shape=circle,draw=black,scale=0.5] (C) at (1.8,0) {A};
\node[shape=circle,draw=black,scale=0.5] (D) at (0,0) {A};

\path [->](A) edge node[above] {a} (B);
\path [->](B) edge node[right] {a} (C);
\path [->](C) edge node[below] {b} (D);
\path [->](D) edge node[left] {b} (A);
\end{tikzpicture}
\caption{$a a b b $}
\end{subfigure}
 \hfill 
\begin{subfigure}{.18\linewidth}
\centering
\begin{tikzpicture}[>={Latex[width=1mm,length=1.5mm]},
                               node distance = 1.5cm and 2cm,
                               el/.style = {inner sep=2pt, align=center, sloped},
                               every label/.append style = {font=\tiny}]
\node[shape=circle,draw=black,scale=0.5] (A) at (0,1.8) {A};
\node[shape=circle,draw=black,scale=0.5] (B) at (1.8,1.8) {B};
\node[shape=circle,draw=black,scale=0.5] (C) at (1.8,0) {B};
\node[shape=circle,draw=black,scale=0.5] (D) at (0,0) {A};

\path [->](A) edge node[above] {a} (B);
\path [->](B) edge node[right] {b} (C);
\path [<-](C) edge node[below] {a} (D);
\path [->](D) edge node[left] {c} (A);
\end{tikzpicture}
\caption{$a b a^{-1} c $}
\end{subfigure}
 \hfill 
\begin{subfigure}{.18\linewidth}
\centering
\begin{tikzpicture}[>={Latex[width=1mm,length=1.5mm]},
                               node distance = 1.5cm and 2cm,
                               el/.style = {inner sep=2pt, align=center, sloped},
                               every label/.append style = {font=\tiny}]
\node[shape=circle,draw=black,scale=0.5] (A) at (0,1.8) {A};
\node[shape=circle,draw=black,scale=0.5] (B) at (1.8,1.8) {A};
\node[shape=circle,draw=black,scale=0.5] (C) at (1.8,0) {A};
\node[shape=circle,draw=black,scale=0.5] (D) at (0,0) {A};

\path [->](A) edge node[above] {a} (B);
\path [->](B) edge node[right] {b} (C);
\path [<-](C) edge node[below] {a} (D);
\path [<-](D) edge node[left] {b} (A);
\end{tikzpicture}
\caption{$a b a^{-1} b^{-1} $}
\end{subfigure}
 \hfill 
\begin{subfigure}{.18\linewidth}
\centering
\begin{tikzpicture}[>={Latex[width=1mm,length=1.5mm]},
                               node distance = 1.5cm and 2cm,
                               el/.style = {inner sep=2pt, align=center, sloped},
                               every label/.append style = {font=\tiny}]
\node[shape=circle,draw=black,scale=0.5] (A) at (0,1.8) {A};
\node[shape=circle,draw=black,scale=0.5] (B) at (1.8,1.8) {A};
\node[shape=circle,draw=black,scale=0.5] (C) at (1.8,0) {A};
\node[shape=circle,draw=black,scale=0.5] (D) at (0,0) {A};

\path [->](A) edge node[above] {a} (B);
\path [->](B) edge node[right] {b} (C);
\path [<-](C) edge node[below] {a} (D);
\path [->](D) edge node[left] {b} (A);
\end{tikzpicture}
\caption{$a b a^{-1} b $}
\end{subfigure}
 \hfill 
\begin{subfigure}{.18\linewidth}
\centering
\begin{tikzpicture}[>={Latex[width=1mm,length=1.5mm]},
                               node distance = 1.5cm and 2cm,
                               el/.style = {inner sep=2pt, align=center, sloped},
                               every label/.append style = {font=\tiny}]
\node[shape=circle,draw=black,scale=0.5] (A) at (0,1.8) {A};
\node[shape=circle,draw=black,scale=0.5] (B) at (1.8,1.8) {B};
\node[shape=circle,draw=black,scale=0.5] (C) at (1.8,0) {A};
\node[shape=circle,draw=black,scale=0.5] (D) at (0,0) {B};

\path [->](A) edge node[above] {a} (B);
\path [->](B) edge node[right] {b} (C);
\path [->](C) edge node[below] {a} (D);
\path [->](D) edge node[left] {c} (A);
\end{tikzpicture}
\caption{$a b a c $}
\end{subfigure}
 \hfill 
\begin{subfigure}{.18\linewidth}
\centering
\begin{tikzpicture}[>={Latex[width=1mm,length=1.5mm]},
                               node distance = 1.5cm and 2cm,
                               el/.style = {inner sep=2pt, align=center, sloped},
                               every label/.append style = {font=\tiny}]
\node[shape=circle,draw=black,scale=0.5] (A) at (0,1.8) {A};
\node[shape=circle,draw=black,scale=0.5] (B) at (1.8,1.8) {B};
\node[shape=circle,draw=black,scale=0.5] (C) at (1.8,0) {A};
\node[shape=circle,draw=black,scale=0.5] (D) at (0,0) {B};

\path [->](A) edge node[above] {a} (B);
\path [->](B) edge node[right] {b} (C);
\path [->](C) edge node[below] {a} (D);
\path [->](D) edge node[left] {b} (A);
\end{tikzpicture}
\caption{$a b a b $}
\end{subfigure}
 \hfill 

    \caption{All quadrilateral gluings}
    \label{fig:quads}
\end{figure}

For $n=1$ the problem is trivial and all gluings can be easily enumerated \ref{fig:quads}. However, that is not the case for $n>1$.

\par
\begin{figure}[!ht]
    \begin{subfigure}{.2\linewidth}
\centering
\begin{tikzpicture}[>={Latex[width=1mm,length=1.5mm]},
                               node distance = 1.5cm and 2cm,
                               el/.style = {inner sep=1pt, align=center, sloped},
                               every label/.append style = {font=\tiny}]
\node[shape=circle,draw=black, scale=0.5] (A) at (0,2) {A};
\node[shape=circle,draw=black, scale=0.5] (B) at (1,2) {B};
\node[shape=circle,draw=black, scale=0.5] (C) at (1.75,1) {C};
\node[shape=circle,draw=black, scale=0.5] (D) at (1,0) {D};
\node[shape=circle,draw=black, scale=0.5] (E) at (0,0) {E};
\node[shape=circle,draw=black, scale=0.5] (F) at (-0.75,1) {F};

\path [->](A) edge node[above] {a} (B);
\path [->](B) edge node[above right] {b} (C);
\path [->](C) edge node[below right] {c} (D);
\path [->](D) edge node[below] {d} (E);
\path [->](E) edge node[below left] {e} (F);
\path [->](F) edge node[above left] {f} (A);
\path [-](A) edge[color=black] node {} (D);
\end{tikzpicture}
\caption{Identity}
\end{subfigure}%
    \hfill
\begin{subfigure}{.2\linewidth}
\centering
\begin{tikzpicture}[>={Latex[width=1mm,length=1.5mm]},
                               node distance = 1.5cm and 2cm,
                               el/.style = {inner sep=1pt, align=center, sloped},
                               every label/.append style = {font=\tiny}]
\node[shape=circle,draw=black, scale=0.5] (A) at (0,2) {D};
\node[shape=circle,draw=black, scale=0.5] (B) at (1,2) {E};
\node[shape=circle,draw=black, scale=0.5] (C) at (1.75,1) {F};
\node[shape=circle,draw=black, scale=0.5] (D) at (1,0) {A};
\node[shape=circle,draw=black, scale=0.5] (E) at (0,0) {B};
\node[shape=circle,draw=black, scale=0.5] (F) at (-0.75,1) {C};

\path [->](A) edge node[above] {d} (B);
\path [->](B) edge node[above right] {e} (C);
\path [->](C) edge node[below right] {f} (D);
\path [->](D) edge node[below] {a} (E);
\path [->](E) edge node[below left] {b} (F);
\path [->](F) edge node[above left] {c} (A);
\path [-](A) edge[color=black] node {} (D);
\end{tikzpicture}
\caption{$\pi$ rotation}
\end{subfigure}
    \hfill
\begin{subfigure}{.2\linewidth}
\centering
\begin{tikzpicture}[>={Latex[width=1mm,length=1.5mm]},
                               node distance = 1.5cm and 2cm,
                               el/.style = {inner sep=1pt, align=center, sloped},
                               every label/.append style = {font=\tiny}]
\node[shape=circle,draw=black, scale=0.5] (A) at (0,2) {A};
\node[shape=circle,draw=black, scale=0.5] (B) at (1,2) {F};
\node[shape=circle,draw=black, scale=0.5] (C) at (1.75,1) {E};
\node[shape=circle,draw=black, scale=0.5] (D) at (1,0) {D};
\node[shape=circle,draw=black, scale=0.5] (E) at (0,0) {C};
\node[shape=circle,draw=black, scale=0.5] (F) at (-0.75,1) {B};

\path [<-](A) edge node[above] {f} (B);
\path [<-](B) edge node[above right] {e} (C);
\path [<-](C) edge node[below right] {d} (D);
\path [<-](D) edge node[below] {c} (E);
\path [<-](E) edge node[below left] {b} (F);
\path [<-](F) edge node[above left] {a} (A);
\path [-](A) edge[color=black] node {} (D);
\end{tikzpicture}
\caption{reflection}
\end{subfigure}
\hfill
\begin{subfigure}{.2\linewidth}
\centering
\begin{tikzpicture}[>={Latex[width=1mm,length=1.5mm]},
                               node distance = 1.5cm and 2cm,
                               el/.style = {inner sep=1pt, align=center, sloped},
                               every label/.append style = {font=\tiny}]
\node[shape=circle,draw=black, scale=0.5] (A) at (0,2) {D};
\node[shape=circle,draw=black, scale=0.5] (B) at (1,2) {C};
\node[shape=circle,draw=black, scale=0.5] (C) at (1.75,1) {B};
\node[shape=circle,draw=black, scale=0.5] (D) at (1,0) {A};
\node[shape=circle,draw=black, scale=0.5] (E) at (0,0) {F};
\node[shape=circle,draw=black, scale=0.5] (F) at (-0.75,1) {E};

\path [<-](A) edge node[above] {c} (B);
\path [<-](B) edge node[above right] {b} (C);
\path [<-](C) edge node[below right] {a} (D);
\path [<-](D) edge node[below] {f} (E);
\path [<-](E) edge node[below left] {e} (F);
\path [<-](F) edge node[above left] {d} (A);
\path [-](A) edge[color=black] node {} (D);
\end{tikzpicture}
\caption{reflection and $\pi$ rotation}
\end{subfigure}
    \caption{Hexagon gluing configuration symmetries}
    \label{fig:hex_conf}
\end{figure}
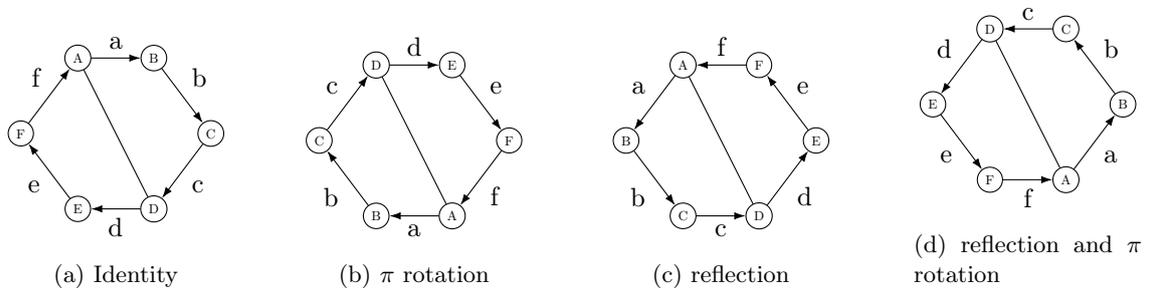

For $n=2$ we have only one possible configuration of hexagon glued from two quadrilaterals. Additionally, this configuration has 6 possible pivot points and 4 symmetries including identity \ref{fig:hex_conf}. Changing a pivot point for gluing configuration correspond to a rotation of a polygon.

\par
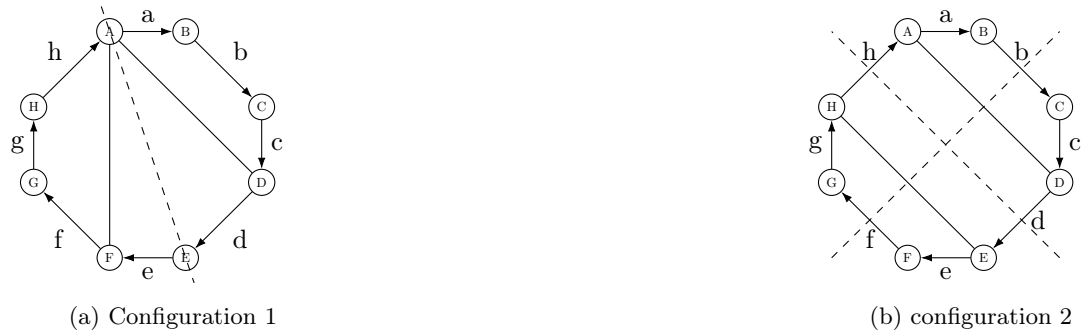
\begin{figure}[!ht]
\begin{subfigure}{.3\linewidth}
\begin{tikzpicture}[>={Latex[width=1.000mm,length=1.500mm]},
                                node distance = 1.500cm and 2.000cm,
                                el/.style = {inner sep=2pt, align=center, sloped},
                                every label/.append style = {font=\tiny}]
\node[shape=circle,draw=black,scale=0.500] (A) at (-0.000,3.000) {A};
\node[shape=circle,draw=black,scale=0.500] (B) at (1.000,3.000) {B};
\node[shape=circle,draw=black,scale=0.500] (C) at (2.000,2.000) {C};
\node[shape=circle,draw=black,scale=0.500] (D) at (2.000,1.000) {D};
\node[shape=circle,draw=black,scale=0.500] (E) at (1.000,0.000) {E};
\node[shape=circle,draw=black,scale=0.500] (F) at (0.000,0.000) {F};
\node[shape=circle,draw=black,scale=0.500] (G) at (-1.000,1.000) {G};
\node[shape=circle,draw=black,scale=0.500] (H) at (-1.000,2.000) {H};

\path [->](A) edge node[above] {a} (B);
\path [->](B) edge node[above right] {b} (C);
\path [->](C) edge node[right] {c} (D);
\path [->](D) edge node[below right] {d} (E);
\path [->](E) edge node[below] {e} (F);
\path [->](F) edge node[below left] {f} (G);
\path [->](G) edge node[left] {g} (H);
\path [->](H) edge node[above left] {h} (A);
\path [-](A) edge[color=black] node {} (D);
\path [-](A) edge[color=black] node {} (F);
\path [dashed](-0.111, 3.333) edge node {} (1.111, -0.333);
\end{tikzpicture}
\caption{Configuration 1}
\label{fig:oct_conf_1}
\end{subfigure}
\hfill
\begin{subfigure}{.3\linewidth}
\begin{tikzpicture}[>={Latex[width=1.000mm,length=1.500mm]},
                                node distance = 1.500cm and 2.000cm,
                                el/.style = {inner sep=2pt, align=center, sloped},
                                every label/.append style = {font=\tiny}]
\node[shape=circle,draw=black,scale=0.500] (A) at (0.000,3.000) {A};
\node[shape=circle,draw=black,scale=0.500] (B) at (1.000,3.000) {B};
\node[shape=circle,draw=black,scale=0.500] (C) at (2.000,2.000) {C};
\node[shape=circle,draw=black,scale=0.500] (D) at (2.000,1.000) {D};
\node[shape=circle,draw=black,scale=0.500] (E) at (1.000,0.000) {E};
\node[shape=circle,draw=black,scale=0.500] (F) at (0.000,0.000) {F};
\node[shape=circle,draw=black,scale=0.500] (G) at (-1.000,1.000) {G};
\node[shape=circle,draw=black,scale=0.500] (H) at (-1.000,2.000) {H};

\path [->](A) edge node[above] {a} (B);
\path [->](B) edge node[above] {b} (C);
\path [->](C) edge node[right] {c} (D);
\path [->](D) edge node[right] {d} (E);
\path [->](E) edge node[below] {e} (F);
\path [->](F) edge node[below] {f} (G);
\path [->](G) edge node[left] {g} (H);
\path [->](H) edge node[above] {h} (A);

\path [-](A) edge[color=black] node {} (D);
\path [-](E) edge[color=black] node {} (H);
\draw[dashed] (-1,0) to (2, 3);
\draw[dashed] (2,0) to (-1, 3);
\end{tikzpicture}
\caption{configuration 2}
\label{fig:oct_conf_2}
\end{subfigure}
    \caption{Octagon gluing configurations}
\end{figure}

For $n=3$ there are two possible configurations of octagon glued from three quadrilaterals \ref{fig:oct_conf_1} and \ref{fig:oct_conf_2} each consisting of two edges. Configuration \ref{fig:oct_conf_2} has 4 symmetries and configuration \ref{fig:oct_conf_1} has only two, and each configuration has 8 possible pivot points.

\par
We consider two gluings distinct if the edge is between different vertices and if labeling schemes are distinct up to relabeling under all symmetries mapping configuration to the same pivot point.

\section{Algorithms}
For small $n$ the problem can be solved using non-optimized but easy-to-implement recursive brute-force algorithms, here we give a brief description of algorithms used to obtain the results for $n=1,2,3$.
\par
First of all, we need to have standard labeling to check whether two labeling schemes are equal up to relabeling.

\begin{algorithm}[!ht]
\caption{Relabeling}\label{relabel}
\begin{algorithmic}[1]
    \Procedure{relabel}{$w$}\Comment{Relabel scheme $w$ to it's standard labeling $w'$}
    \State $w' \gets empty(w)$\Comment{Empty scheme of length $w$}
    \For{i:=1 \textbf{to} $length(w)$}
        \If{$isset(w(i))$}\Comment{Letter at $i$ is present}
          \State $continue$
        \Else
          \State $e \gets nextLetter(w')$\Comment{Next alphabet letter that is not in $w'$}
          \State $w'(i) \gets e$
          \If{$occurrences(w,w(i)) > 1$}
            \State $j \gets gluedIndex(w,i)$
            \If{$exponent(w(i))=exponent(w(j))$}
            \State $w'(j) \gets e$
            \Else
            \State $w'(j) \gets inverse(w(e))$
            \EndIf
          \EndIf
        \EndIf
      \EndFor
    \State \textbf{return} $w'$\Comment{Standard labeling of $w$}
\EndProcedure
\end{algorithmic}
\end{algorithm}

Examples of relabeling algorithm inputs and outputs:
$$b a^{-1} c d c^{-1} d^{-1} b^{-1} a \mapsto a b c d c^{-1} d^{-1} a^{-1} b^{-1}$$
$$a^{-1} b^{-1} a d^{-1} e^{-1} f g^{-1} f \mapsto a b a^{-1} c d e f e  $$

\begin{algorithm}[!ht]
\caption{Edge gluing}\label{gluing}
\begin{algorithmic}[1]
    \Procedure{glue}{$w,i_1,i_2,orientable?$}\Comment{Glue edges at indices $i_1$ and $i_2$ of labeling scheme $w$}
    \State $e \gets w(i_1)$
    \If {$orientable? = true$} 
        \State $w(i_2) \gets inverse(e)$
    \Else
        \State $w(i_2) \gets e$
    \EndIf
    \State $w \gets relabel(w)$
    \State \textbf{return} $w$\Comment{Labeling scheme with edges at $i_1$ and $i_2$ glued}
\EndProcedure
\end{algorithmic}
\end{algorithm}

The vertices labeling algorithm is a crucial part of finding the Euler characteristic of a surface defined by the labeling scheme\ref{vertices}

\begin{algorithm}[ht!]
\caption{Vertices labeling}\label{vertices}
\begin{algorithmic}[1]
    \Procedure{vertices}{$w$}\Comment{Get vertices sequence for $w$}
    \State $V \gets empty(w)$\Comment{Empty collection of vertices same size as $w$}
    \Function{$updateEquivalence$}{$eq,w,i$}
    \For{$e$ : $adjacentEdges(w, i)$}
        \If{$isGlued(e)$}
        \State $j \gets gluedVertexIdx(e,i)$
        \State $addVertex(eq, j)$
        \State $eq \gets updateEquivalence(eq,w,j)$\Comment{Recursive update for glued vertex}
        \EndIf
    \EndFor
    \State \textbf{return} $eq$
    \EndFunction
    \For{i:=1 \textbf{to} $length(w)$}
    \If{$V(i)$}
    \State $continue$
    \EndIf
    \State $eq := set()$\Comment{Equivalence class of vertex indices}
    \State $eq \gets updateEquivalence(eq,w,i)$
    \State $V \gets setVertices(V, eq)$\Comment{Set indices from equivalence class to next alphabet letter not in $V$}
    \EndFor
    \State \textbf{return} $V$
\EndProcedure
\end{algorithmic}
\end{algorithm}
\newpage

\section{Classifications of glueings}
For $n = 1$ we have 7 topologically different surfaces:
\begin{table}[!ht]
\begin{tabular}{|c|c|c|c|c|c|c|c|}
\hline
$\chi$              & 2  & 0    & 1    & 0    & 1     & 0     & 0     \\
\hline
orientable       & true  & true & true & true & false & false & false \\
\hline
boundary  & 0  & 0    & 1    & 2    & 0     & 0     & 1      \\
\hline
n                   & 1  & 1    & 1    & 1    & 2     & 2     & 2       \\
\hline
\end{tabular}
\end{table}

\par
For $n = 2$ we have 11 topologically different structures, the number of different graphs embedded on each surface is calculated with respect to different configurations of gluing fundamental hexagon from two quadrilaterals.

\begin{table}[!ht]
\begin{tabular}{|c|c|c|c|c|c|c|c|c|c|c|c|}
\hline
$\chi$ & 0 & 1 & 1 & -1 & 2 & -1 & -1 & 0 & 0 & -1 & 0\\ 
 \hline 
orientable & true & true & false & false & true & true & false & false & false & false & true\\ 
 \hline 
boundary & 0 & 1 & 0 & 1 & 0 & 1 & 2 & 1 & 0 & 0 & 2\\ 
 \hline 
n & 5 & 9 & 10 & 15 & 3 & 5 & 7 & 22 & 14 & 8 & 10\\ 
 \hline 
\end{tabular}
\end{table}
\par

For $n = 3$ we have many more possible gluing configurations and 17 non-homeomorphic surfaces.

\begin{table}[!ht]
\begin{tabular}{|c|c|c|c|c|c|c|c|c|c|c|c|c|c|c|c|c|c|}
\hline
$\chi$ & -1 & 1 & 0 & -1 & 1 & 0 & 0 & -1 & -2 & -1 & -2 & 2 & -1 & -2 & 0 & -2 & -2\\ 
 \hline 
o & t & f & f & f & t & t & t & f & f & t & f & t & f & f & f & t & t\\ 
 \hline 
b & 3 & 0 & 0 & 2 & 1 & 2 & 0 & 0 & 0 & 1 & 2 & 0 & 1 & 1 & 1 & 0 & 2\\ 
 \hline 
n & 19 & 59 & 145 & 219 & 65 & 113 & 43 & 197 & 47 & 122 & 95 & 16 & 451 & 256 & 263 & 4 & 25\\ 
 \hline 
\end{tabular}
\end{table}

\section{Conclusion}
In this paper, we have developed a method for classification graphs consisting of bounding quadrilateral regions on Morse fields (Morse foliations). The result of this work can be applied in many fields. For example, these kinds of foliations, are created by Morse flows.


\begin{thebibliography}{10}

\bibitem{akchurin2022three}
O.~Akchurin, S.~Bilun, and A.~Prishlyak.
\newblock Three-color graph as the 1-skeleton of the 2-sphere triangulation.
\newblock {\em arXiv preprint arXiv:2209.05737}, 2022.
\href{http://dx.doi.org/10.48550/ARXIV.2209.05737 }{\path{doi: 10.48550/ARXIV.2209.05737}}.

\bibitem{bilun2002closed}
S.~Bilun and A.~Prishlyak.
\newblock The closed morse 1-forms on closed surfaces.
\newblock {\em Visn., Mat. Mekh., Kyv. Univ. Im. Tarasa Shevchenka},
  2002(8):77--81, 2002.

\bibitem{bilun2022visualization}
S.~Bilun and A.~Prishlyak.
\newblock Visualization of morse flow with two saddles on 3-sphere diagrams.
\newblock {\em arXiv preprint arXiv:2209.12174}, 2022.
\href{http://dx.doi.org/10.48550/ARXIV.2209.12174}{\path{doi: 10.48550/ARXIV.2209.12174}}.


\bibitem{bilun2022morse}
S.~Bilun, A.~Prishlyak, and A.~Prus.
\newblock Morse flows with fixed points on the boundary of 3-manifold.
\newblock {\em arXiv preprint arXiv:2209.04019}, 2022.
\href{http://dx.doi.org/10.48550/arXiv.2209.04019 }{\path{doi: 10.48550/arXiv.2209.04019}}.


\bibitem{bilun2023discrete}
S. Bilun, M. Hrechko, O. Myshnova, A. Prishlyak.
\newblock Structures of optimal discrete gradient vector fields on surface with one or two critical cells
\newblock {\em arXiv preprint arXiv:2303.07258}, 2023.
\href{http://dx.doi.org/10.48550/arXiv.2303.07258}{\path{doi: 10.48550/arXiv.2303.07258}}.


\bibitem{bilun2023morseRP2}
S. Bilun, A. Prishlyak, S. Stas, A. Vlasenko.
\newblock Topological structure of Morse functions on the projective plane
\newblock {\em arXiv preprint arXiv:2303.03850}, 2023.
\href{http://dx.doi.org/10.48550/arXiv.2303.03850}{\path{doi: 10.48550/arXiv.2303.03850}}.


\bibitem{bilun2023gradient}
S. Bilun, B. Hladysh, A. Prishlyak, V. Sinitsyn.
\newblock Gradient vector fields of codimension one on the 2-sphere with at most ten singular points
\newblock {\em arXiv preprint arXiv:2303.10929}, 2023.
\href{http://dx.doi.org/10.48550/arXiv.2303.10929}{\path{doi: 10.48550/arXiv.2303.10929}}.


\bibitem{bilun2023typical}
S. Bilun, M. Loseva, O. Myshnova, A. Prishlyak.
\newblock Typical one-parameter bifurcations of gradient flows with at most six singular points on the 2-sphere with holes
\newblock {\em arXiv preprint arXiv:2303.14975}, 2023.
\href{http://dx.doi.org/10.48550/arXiv.2303.14975}{\path{doi: 10.48550/arXiv.2303.14975}}.


\bibitem{bilun2023morse}
S. Bilun, B. Hladysh, A. Prishlyak, M. Roman.
\newblock Morse functions with four critical points on immersed 2-spheres.
\newblock {\em arXiv preprint arXiv:2304.04392}, 2023.
\href{http://dx.doi.org/10.48550/arXiv.2304.04392}{\path{doi: 10.48550/arXiv.2304.04392}}.


\bibitem{Bolsinov2004}
A.V. Bolsinov and A.T. Fomenko.
\newblock {\em Integrable Hamiltonian systems. Geometry, Topology,
  Classification}.
\newblock A CRC Press Company, Boca Raton London New York Washington, D.C.,
  2004.
\newblock 724 p.


\bibitem{borodzik2016morse}
M. Borodzik, A. Némethi and A. Ranicki
\newblock Morse theory for manifolds with boundary.
\newblock {\em Algebraic Geometric Topology}, 16(2):971--1023, 2016.
\href{http://dx.doi.org/ 10.2140/agt.2016.16.971}{\path{doi:  10.2140/agt.2016.16.971}}.




\bibitem{GM88}
M. Goresky and R. MacPherson.
\newblock {\em Stratified Morse Theory}.
\newblock Springer, 1988.
\href{http://dx.doi.org/10.1007/978-3-642-71714-7}{\path{doi: 10.1007/978-3-642-71714-7}}.



\bibitem{hatamian2020heegaard}
C.~Hatamian and A.~Prishlyak.
\newblock Heegaard diagrams and optimal morse flows on non-orientable
  3-manifolds of genus 1 and genus 2.
\newblock {\em Proceedings of the International Geometry Center}, 13(3):33--48,
  2020.
\href{http://dx.doi.org/10.15673/tmgc.v13i3.1779}{\path{doi: 10.15673/tmgc.v13i3.1779}}.


\bibitem{hladysh2016functions}
B.~I. Hladysh and A.~O. Pryshlyak.
\newblock Functions with nondegenerate critical points on the boundary of the
  surface.
\newblock {\em Ukrainian Mathematical Journal}, 68(1):29--41, 2016.
\href{http://dx.doi.org/10.1007/s11253-016-1206-5}{\path{doi: 10.1007/s11253-016-1206-5}}.

\bibitem{hladysh2020deformations}
B.~I. Hladysh and A.~O. Pryshlyak.
\newblock Deformations in the general position of the optimal functions on oriented surfaces with boundary.
\newblock {\em Ukrainian Mathematical Journal}, 71(8):1173--1185, 2020.
\href{http://dx.doi.org/10.1007/s11253-019-01706-8}{\path{doi: 10.1007/s11253-019-01706-8}}.


\bibitem{hladysh2017topology}
B.I. Hladysh and A.O. Prishlyak.
\newblock Topology of functions with isolated critical points on the boundary
  of a 2-dimensional manifold.
\newblock {\em SIGMA. Symmetry, Integrability and Geometry: Methods and
  Applications}, 13:050, 2017.
\href{http://dx.doi.org/0.3842/SIGMA.2017.050}{\path{doi: 0.3842/SIGMA.2017.050}}.

\bibitem{hladysh2019simple}
B.I. Hladysh and A.O. Prishlyak.
\newblock Simple morse functions on an oriented surface with boundary.
\newblock {\em Журнал математической физики,
  анализа, геометрии}, 15(3):354--368, 2019.
\href{http://dx.doi.org/10.15407/mag15.03.354}{\path{doi: 10.15407/mag15.03.354}}.

\bibitem{Kronrod1950}
A.S. Kronrod.
\newblock Functions of two variables.
\newblock {\em Russian Mathematical Surveys}, 5:24--134, 1950.


\bibitem{Kybalko2018}
Z.~Kybalko, A.~Prishlyak, and R.~Shchurko.
\newblock {Trajectory equivalence of optimal Morse flows on closed surfaces}.
\newblock {\em {Proc. Int. Geom. Cent.}}, 11(1):12--26, 2018.
\href{http://dx.doi.org/10.15673/tmgc.v11i1.916 }{\path{doi: 10.15673/tmgc.v11i1.916}}.

\bibitem{loseva2022topological}
M.~Loseva, A.~Prishlyak.
\newblock Topological structure of optimal flows on the Girl’s surface.
\newblock {\em Proceedings of the International Geometry Center}, 15(3-4):184--202,
  2022.
\href{https://doi.org/10.15673/tmgc.v15i3-4.2338}{\path{doi: 10.15673/tmgc.v15i3-4.2338}}.

\bibitem{loseva2016topology}
M.~Losieva and A.~Prishlyak.
\newblock Topology of morse--smale flows with singularities on the boundary of
  a two-dimensional disk.
\newblock {\em Pr. Mizhnar. Heometr. Tsentr}, 9(2):32--41, 2016.
\href{http://dx.doi.org/10.15673/tmgc.v9i2.279}{\path{doi: 10.15673/tmgc.v9i2.279}}.




\bibitem{lychak2009morse}
D.P. Lychak and A.O. Prishlyak.
\newblock Morse functions and flows on nonorientable surfaces.
\newblock {\em Methods of Functional Analysis and Topology}, 15(03):251--258,
  2009.



\bibitem{Oshemkov1998}
A.A. Oshemkov and V.V. Sharko.
\newblock Classication of morse-smale flows on two-dimensional manifolds.
\newblock {\em Matem. Sbornik}, 189(8):93--140, 1998.

\bibitem{Peixoto1973}
M.M. Peixoto.
\newblock On the classication of flows of 2-manifolds.
\newblock {\em Dynamical Systems (Proc. Symp. Univ. of Bahia, Salvador, Brasil,
  1971)}, pages 389--419, 1973.


\bibitem{prishlyak1997graphs}
A.O. Prishlyak.
\newblock On graphs embedded in a surface.
\newblock {\em Russian Mathematical Surveys}, 52(4):844, 1997.
\href{http://dx.doi.org/10.1070/RM1997v052n04ABEH002074}{\path{doi: 10.1070/RM1997v052n04ABEH002074}}.

\bibitem{prishlyak2001conjugacy}
A.O. Prishlyak.
\newblock Conjugacy of Morse functions on 4-manifolds.
\newblock {\em Russian Mathematical Surveys}, 56(1):170, 2001.




\bibitem{prishlyak2002ms}
A.O. Prishlyak.
\newblock Morse--smale vector fields without closed trajectories on-manifolds.
\newblock {\em Mathematical Notes}, 71(1-2):230--235, 2002.
\href{http://dx.doi.org/10.1023/A:1013963315626}{\path{doi: 10.1023/A:1013963315626}}.

\bibitem{prishlyak2003sum}
A.O. Prishlyak.
\newblock On sum of indices of flow with isolated fixed points on a stratified
  sets.
\newblock {\em Zhurnal Matematicheskoi Fiziki, Analiza, Geometrii [Journal of
  Mathematical Physics, Analysis, Geometry]}, 10(1):106--115, 2003.


\bibitem{prishlyak2007complete}
A.O. Prishlyak.
\newblock Complete topological invariants of morse-smale flows and handle
  decompositions of 3-manifolds.
\newblock {\em Journal of Mathematical Sciences}, 144:4492--4499, 2007.

\bibitem{prishlyak2012topological}
A.O. Prishlyak.
\newblock Topological properties of functions on two and three dimensional manifolds.
\newblock {\em {Palmarium Academic Pablishing}}, 2012.


\bibitem{prishlyak2022Boy}
A.~Prishlyak and L.~Di Beo.
\newblock Flows with minimal number of singularities on the Boy’s surface
\newblock {\em Proceedings of the International Geometry Center}, 15(1):32--49,
  2020.

\bibitem{prishlyak2022topological}
A.~Prishlyak and M.~Loseva.
\newblock Topological structure of optimal flows on the girl's surface.
\newblock {\em Proceedings of the International Geometry Center},
  15(3-4):184--202, 2022.

\bibitem{prishlyak2021flows}
A.~Prishlyak, A.~Prus, and S.~Huraka.
\newblock Flows with collective dynamics on a sphere.
\newblock {\em Proc. Int. Geom. Cent}, 14(1):61--80, 2021.
\href{http://dx.doi.org/10.15673/tmgc.v14i1.1902}{\path{doi: 10.15673/tmgc.v14i1.1902}}.




\bibitem{prishlyak2020topology}
A.~Prishlyak and M.~Loseva.
\newblock Topology of optimal flows with collective dynamics on closed
  orientable surfaces.
\newblock {\em Proceedings of the International Geometry Center}, 13(2):50--67,
  2020.
\href{http://dx.doi.org/10.15673/tmgc.v13i2.1731}{\path{doi: 10.15673/tmgc.v13i2.1731}}.

\bibitem{prishlyak2017morse}
A.~Prishlyak and A.~Prus.
\newblock Morse-smale flows on torus with hole.
\newblock {\em Proc. Int. Geom. Cent.}, 10(1):47--58, 2017.
\href{http://dx.doi.org/10.15673/tmgc.v1i10.549}{\path{doi: 10.15673/tmgc.v1i10.549}}.

\bibitem{prishlyak1999equivalence}
A.O.~Prishlyak.
\newblock Equivalence of morse function on 3-manifolds.
\newblock {\em Methods of Func. Ann. and Topology}, 5(3):49--53, 1999.

\bibitem{prishlyak2000conjugacy}
A.O.~Prishlyak.
\newblock Conjugacy of morse functions on surfaces with values on a straight
  line and circle.
\newblock {\em Ukrainian Mathematical Journal}, 52(10):1623--1627, 2000.
\href{http://dx.doi.org/10.1023/A:1010461319703}{\path{doi: 10.1023/A:1010461319703}}.

\bibitem{prishlyak2002morse}
A.O.~Prishlyak.
\newblock Morse functions with finite number of singularities on a plane.
\newblock {\em Meth. Funct. Anal. Topol}, 8:75--78, 2002.

\bibitem{Prishlyak2002beh2}
A.O.~Prishlyak.
\newblock Topological equivalence of morse--smale vector fields with beh2 on
  three-dimensional manifolds.
\newblock {\em Ukrainian Mathematical Journal}, 54(4):603--612, 2002.

\bibitem{prishlyak2003topological}
A.O.~Prishlyak.
\newblock Topological classification of m-fields on two-and three-dimensional
  manifolds with boundary.
\newblock {\em Ukrainian Mathematical Journal}, 55(6):966--973, 2003.

\bibitem{prishlyak2019optimal}
A.O.~Prishlyak and M.V.~Loseva.
\newblock Optimal morse--smale flows with singularities on the boundary of a
  surface.
\newblock {\em Journal of Mathematical Sciences}, 243:279--286, 2019.

\bibitem{prishlyak2007classification}
A.O.~Prishlyak and K.I.~Mischenko.
\newblock Classification of noncompact surfaces with boundary.
\newblock {\em Methods of Functional Analysis and Topology}, 13(01):62--66,
  2007.

\bibitem{prishlyak2020three}
A.O.~Prishlyak and A.A.~Prus.
\newblock Three-color graph of the morse flow on a compact surface with
  boundary.
\newblock {\em Journal of Mathematical Sciences}, 249(4):661--672, 2020.
\href{http://dx.doi.org/10.1007/s10958-020-04964-1}{\path{doi: 10.1007/s10958-020-04964-1}}.

\bibitem{stas2023structures}
A.~Prishlyak and S.~Stas.
\newblock Structures of the flows with a unique singular point on the 2-dimensional disk.
\newblock {\em arXiv preprint arXiv:2304.00751}, 2023.
\href{http://dx.doi.org/10.48550/ARXIV.2304.00751}{\path{doi: 10.48550/ARXIV.2304.00751}}.



\bibitem{Reeb1946}
G.~Reeb.
\newblock Sur les points singuliers d’une forme de pfaff complétement
  intégrable ou d’une fonction numérique.
\newblock {\em C.R.A.S. Paris}, 222:847—849, 1946.

\bibitem{Sharko1993}
V.V. Sharko.
\newblock {\em Functions on manifolds. Algebraic and topological aspects.},
  volume 131 of {\em Translations of Mathematical Monographs}.
\newblock American Mathematical Society, Providence, RI, 1993.

\bibitem{Smale1961}
S.~Smale.
\newblock On gradient dynamical systems.
\newblock {\em Ann. of Math.}, 74:199--206, 1961.

\bibitem{kkp2013}
В.М. Кузаконь, В.Ф. Кириченко, and О.О. Пришляк.
\newblock Гладкi многовиди. Геометричнi та
  топологiчнi аспекти.
\newblock {\em Працi Iнcтитуту математики НАН
  України.—2013.—97.—500 с}, 2013.

\bibitem{prish1998vek}
А.О. Пришляк.
\newblock Векторные поля Морса--Смейла с
  конечным числом особых траекторий на
  трехмерных многообразиях.
\newblock {\em Доповiдi НАН України}, (6):43--47, 1998.

\bibitem{prish1998sopr}
А.О. Пришляк.
\newblock Сопряженность функций Морса.
\newblock {\em Некоторые вопросы совр.
  математики. Институт математики АН
  Украины, Киев}, 1998.

\bibitem{prish2001top}
А.О. Пришляк.
\newblock Топологическая эквивалентность
  функций и векторных полей Морса—Смейла
  на трёхмерных многообразиях.
\newblock {\em Топология и геометрия. Труды
  Украинского мат. конгресса}, pages 29--38, 2001.


\bibitem{prish2002theory}
О.О. Пришляк.
\newblock Теорія Морса.
\newblock {\em Київський університет}, 2002.

\bibitem{prish2004difgeom}
О.О. Пришляк.
\newblock Диференціальна геометрія.
\newblock {\em Київський університет}, 2004.


\bibitem{prish2006osnovy}
О.О. Пришляк.
\newblock Основи сучасної топології.
\newblock {\em Київський університет}, 2006.

\bibitem{prish2015top}
О.О. Пришляк.
\newblock Топологія многовидів.
\newblock {\em Київський університет}, 2015.





\bibitem{angelini2018linear}
P.~Angelini, G.~Da Lozzo, G.~Di Battista, F.~Frati, M.~Patrignani, and I.~Rutter.
\newblock Linear area drawings of graphs.
\newblock {\em Algorithmica}, 80(1):174--197, 2018.
\href{http://dx.doi.org/10.1007/s00453-017-0341-8}{\path{doi: 10.1007/s00453-017-0341-8}}.

\bibitem{appel1976every}
K.~Appel and W.~Haken.
\newblock Every planar map is four colorable. {I}. {D}ischarging.
\newblock {\em Illinois Journal of Mathematics}, 21(3):429--490, 1977.
\href{http://dx.doi.org/10.1215/ijm/1256049011}{\path{doi: 10.1215/ijm/1256049011}}.

\bibitem{bang2008digraphs}
J.~Bang-Jensen and G.~Gutin.
\newblock Digraphs: theory, algorithms and applications.
\newblock {\em Springer}, 2008.
\href{http://dx.doi.org/10.1007/978-1-84628-970-5}{\path{doi: 10.1007/978-1-84628-970-5}}.



\bibitem{bollobas1998modern}
B.~Bollob{'a}s.
\newblock Modern graph theory.
\newblock {\em Graduate Texts in Mathematics}, volume 184. Springer, 1998.
\href{http://dx.doi.org/10.1007/978-1-4612-0619-4}{\path{doi: 10.1007/978-1-4612-0619-4}}.


\bibitem{bondy2008graph}
J.A.~Bondy and U.S.R.~Murty.
\newblock Graph theory.
\newblock {\em Graduate Texts in Mathematics}, volume 244. Springer, 2008.
\href{http://dx.doi.org/10.1007/978-1-84628-970-5}{\path{doi: 10.1007/978-1-84628-970-5}}.

\bibitem{chrobak1994planar}
M.~Chrobak and M.T.~Goodrich.
\newblock Planar orthogonal drawings with optimal area.
\newblock {\em Computational Geometry: Theory and Applications}, 4(3):195--204, 1994.
\href{http://dx.doi.org/10.1016/0925-7721(94)90014-0}{\path{doi: 10.1016/0925-7721(94)90014-0}}.

\bibitem{diestel2005graph}
R.~Diestel.
\newblock Graph theory.
\newblock {\em Graduate Texts in Mathematics}, volume 173. Springer, 2005.
\href{http://dx.doi.org/10.1007/978-3-662-48541-6}{\path{doi: 10.1007/978-3-662-48541-6}}.

\bibitem{dillencourt1993graph}
M.B.~Dillencourt, D.E.~Knuth, and M.~Rosenstiehl.
\newblock Graph layout problems.
\newblock {\em SIAM Review}, 35(3):315--347, 1993.
\href{http://dx.doi.org/10.1137/1035095}{\path{doi: 10.1137/1035095}}.


\bibitem{eppstein1999quasi}
D.~Eppstein.
\newblock Quasi-planar graphs.
\newblock {\em Transactions of the American Mathematical Society}, 351(12):5247--5275, 1999.
\href{http://dx.doi.org/10.1090/S0002-9947-99-02342-5}{\path{doi: 10.1090/S0002-9947-99-02342-5}}.

\bibitem{godsil2001algebraic}
C.~Godsil and G.~Royle.
\newblock Algebraic graph theory.
\newblock {\em Graduate Texts in Mathematics}, volume 207. Springer, 2001.
\href{http://dx.doi.org/10.1007/978-1-4757-3797-1}{\path{doi: 10.1007/978-1-4757-3797-1}}.


\bibitem{gross1999topological}
J.L.~Gross and T.W.~Tucker.
\newblock Topological graph theory.
\newblock In L.W.~Beineke and R.J.~Wilson, editors, {\em Topics in Topological Graph Theory}, volume 128 of {\em Encyclopedia of Mathematics and its Applications}, pages 1--101. Cambridge University Press, 1999.
\href{http://dx.doi.org/10.1017/CBO9781139171329}{\path{doi: 10.1017/CBO9781139171329}}.

\bibitem{liotta2018beyond}
G.~Liotta and H.~Meijer.
\newblock Beyond pseudo-3D: new directions for drawing planar graphs.
\newblock In {\em Proceedings of the 26th International Symposium on Graph Drawing and Network Visualization}, pages 317--331, 2018.
\href{http://dx.doi.org/10.1007/978-3-030-04414-5_26}{\path{doi: 10.1007/978-3-030-04414-5_26}}.

\bibitem{lovasz2001graph}
L.~Lov{'a}sz.
\newblock Graph theory.
\newblock {\em Graph Theory}, volume 2 of {\em Encyclopedia of Mathematics and its Applications}, pages 1--434. Cambridge University Press, 2001.
\href{http://dx.doi.org/10.1017/CBO9780511806291}{\path{doi: 10.1017/CBO9780511806291}}.

\bibitem{mohar1991embeddings}
B.~Mohar.
\newblock Embeddings of graphs into surfaces.
\newblock {\em Journal of Combinatorial Theory, Series B}, 51(2):150--163, 1991.
\href{http://dx.doi.org/10.1016/0095-8956(91)90043-N}{\path{doi: 10.1016/0095-8956(91)90043-N}}.

\bibitem{mohar2001graphs}
B.~Mohar.
\newblock Graphs on surfaces.
\newblock {\em Johns Hopkins University Press}, 2001.


\bibitem{mohar1997isomorphism}
B.~Mohar.
\newblock Isomorphism of graphs of small genus after subdivision.
\newblock {\em Journal of Combinatorial Theory, Series B}, 70(1):16--46, 1997.
\href{http://dx.doi.org/10.1006/jctb.1997.1752}{\path{doi: 10.1006/jctb.1997.1752}}.

\bibitem{mohar2001some}
B.~Mohar.
\newblock Some applications of Laplace eigenvalues of graphs.
\newblock In L.W.~Beineke and R.J.~Wilson, editors, {\em Topics in Topological Graph Theory}, volume 128 of {\em Encyclopedia of Mathematics and its Applications}, pages 174--202. Cambridge University Press, 1999.
\href{http://dx.doi.org/10.1017/CBO9781139171329}{\path{doi: 10.1017/CBO9781139171329}}.

\bibitem{neumann2015topological}
A.~Neumann, D.~Prestel, and A.~Schulz.
\newblock Topological graph theory.
\newblock {\em CRC Press}, 2015.

\bibitem{pach2018planar}
J.~Pach, R.~Radoi{\v{c}}i{'c}, and G.~Toth.
\newblock Planar embeddings of graphs on surfaces.
\newblock {\em Graphs and Combinatorics}, 34(1):151--185, 2018.
\href{http://dx.doi.org/10.1007/s00373-017-1833-3}{\path{doi: 10.1007/s00373-017-1833-3}}.

\bibitem{robertson1997graph}
N.~Robertson, D.P.~Sanders, P.D.~Seymour, and R.~Thomas.
\newblock The four-color theorem.
\newblock {\em Journal of Combinatorial Theory, Series B}, 70(1):2--44, 1997.
\href{http://dx.doi.org/10.1006/jctb.1997.1750}{\path{doi: 10.1006/jctb.1997.1750}}.

\bibitem{thomassen2009graphs}
C.~Thomassen.
\newblock Graphs on surfaces: Dualities, embeddings and characterizations.
\newblock {\em Springer Science  Business Media}, 2009.

\bibitem{thomassen1994planar}
C.~Thomassen.
\newblock Planar acyclic graphs.
\newblock {\em Journal of Combinatorial Theory, Series B}, 62(1):180--187, 1994.
\href{http://dx.doi.org/10.1006/jctb.1994.1078}{\path{doi: 10.1006/jctb.1994.1078}}.

\bibitem{thomassen1986planarity}
C.~Thomassen.
\newblock Planarity and duality of finite and infinite graphs.
\newblock {\em Journal of Combinatorial Theory, Series B}, 41(3):291--319, 1986.
\href{http://dx.doi.org/10.1016/0095-8956(86)90089-9}{\path{doi: 10.1016/0095-8956(86)90089-9}}.

\bibitem{tutte1961four}
W.T.~Tutte.
\newblock On the algebraic theory of graph colorings.
\newblock {\em Journal of Combinatorial Theory}, 11(3):267--275, 1961.
\href{http://dx.doi.org/10.1016/S0021-9800(61)80063-4}{\path{doi: 10.1016/S0021-9800(61)80063-4}}.


\bibitem{west2001introduction}
D.B.~West.
\newblock Introduction to graph theory.
\newblock {\em Prentice Hall}, 2001.

\bibitem{Harer1986euler}
J. Harer, D. Zagier
The Euler Characteristic of the Moduli Space of Curves.
\newblock {\em Invent Math}, 85(1):457--485, 1986.
\href{https://doi.org/10.1007/BF01390325}{\path{doi: 10.1007/BF01390325}}

\bibitem{akhmedov2008gluing}
E. T. Akhmedov and Sh. Shakirov
\newblock Gluing of Surfaces with Polygonal Boundaries.
\newblock {\em arXiv preprint arXiv:0712.2448}, 2008
\href{https://doi.org/10.48550/arXiv.0712.2448}{\path{doi: 10.48550/arXiv.0712.2448}}

\bibitem{Kosniowski1980first}
C. Kosniowski
\newblock A first course in algebraic topology.
\newblock {\em Cambridge University Press}, 1980


\end{thebibliography}
\end{document}